\documentclass[conference]{IEEEtran}
\usepackage{cite}
\usepackage{amsmath,amssymb,amsfonts,amsthm}
\usepackage{enumitem}
\usepackage{mathtools}
\usepackage{mathrsfs}
\usepackage[active]{srcltx}
\usepackage{bm}
\usepackage{fixmath}
\usepackage{graphicx}
\usepackage{textcomp}
\usepackage{xcolor}
\def\BibTeX{{\rm B\kern-.05em{\sc i\kern-.025em b}\kern-.08em
    T\kern-.1667em\lower.7ex\hbox{E}\kern-.125emX}}

\newtheorem{theorem}{Theorem}
\newtheorem{proposition}{Proposition}
\newtheorem{corollary}{Corollary}
\newtheorem*{remark}{Remark}

\begin{document}

%\title{A new Randers space model and its isometric equivalents}

\title{Three isometrically equivalent models of the Finsler-Poincar\'e disk}

\author{\IEEEauthorblockN{\'{A}gnes Mester}
	\IEEEauthorblockA{\textit{Institute of Applied Mathematics} \\
	\textit{\'Obuda University}\\
	Budapest, Hungary \\
	and\\
	\textit{Department of Mathematics and Computer Science}\\
	\textit{Babe\c s-Bolyai University}\\
	Cluj-Napoca, Romania\\
	%mester.agnes@stud.uni-obuda.hu}
	agnes.mester@ubbcluj.ro}	
	\and
	\IEEEauthorblockN{Alexandru Krist\'aly }
	\IEEEauthorblockA{\textit{Institute of Applied Mathematics} \\
	\textit{\'Obuda University}\\
	Budapest, Hungary \\
	and\\
	\textit{Department of Economics}\\
	\textit{ Babe\c s-Bolyai University}\\
	Cluj-Napoca, Romania\\
	kristaly.alexandru@nik.uni-obuda.hu}
%	alex.kristaly@econ.ubbcluj.ro}
}

\maketitle

\begin{abstract}
We present the isometry between the $2$-dimensional Funk model and the Finsler-Poincar\'e disk. Then, we introduce the Finslerian Poincar\'e upper half plane model, which turns out to be also isometrically equivalent to the previous models. 
As application, we state the gapless character of the first eigenvalue for the aforementioned three spaces. 
\end{abstract}


\section{Introduction and main results}
 
The theory of Finsler manifolds can be considered as a generalization of Riemannian geometry, where the Riemannian metric is replaced by a so called Finsler structure, which is induced by a Minkowski norm. Therefore, Finsler geometry provides a natural framework to study anisotropical phenomena, admitting numerous applications in physics and practical problems, see e.g. Antonelli, Ingarden and Matsumoto \cite{AIM94}, Dehkordi\cite{D21}, Gibbons and Warnick \cite{GW11}, Matsumoto \cite{Matsumoto89} and Randers \cite{Randers41}. 

One of the simplest classes of Finsler manifolds are the so called Randers spaces, which have received much attention lately due to Zermelo's famous navigation problem, see Zermelo \cite{Zermelo}. 
More precisely, if $(M, g)$ is a complete
$n$-dimensional ($n \geq 2$) Riemannian manifold, then the Finsler metric $F: TM \rightarrow \mathbb R$ defined as 
\begin{equation*}
F(x,v) = \sqrt{g_x(v,v)} + \beta_x(v), \ \ x \in M , \ v \in T_xM
\end{equation*}
is called a Randers metric whenever $\beta_x$ is a $1$-form on $M$ with $|\beta_x|_g \coloneqq \sqrt{g^*_x(\beta_x,\beta_x)} < 1$
for every $x \in M$, where $g^*$ denotes the co-metric of $g$. As it turns out, every Randers space $(M,F)$ can be obtained as the solution to the Zermelo navigation problem for a suitable choice of $g$ and $\beta_x$, see Bao and Robles \cite{BR}, Bao, Robles and Shen \cite{BRS}, and Shen \cite{Shen03}. Thus every Randers metric can be written as a suitable perturbation
of a Riemannian metric $g$.  

The two typical analytical models of Randers spaces are the following:
\begin{enumerate}
	\item[(\textbf{F}):]\label{P} the Finslerian
	\textit{Funk model} (see Cheng and Shen \cite[Example 2.1.2]{CS12} and Shen \cite[Example 1.3.4]{Shen01}), which turns out to be the generalization of the well known Riemannian \textit{Klein model}; 
	
	\item[(\textbf{P}):]\label{F} the \textit{Finsler-Poincar\'e disk} (see Bao, Chern and Shen \cite[Section 12.6]{BCS}) which appears as the perturbation of the usual Riemannian \textit{Poincar\'e metric} on the open unit disk. 
\end{enumerate} 

As it turns out, these two Randers spaces above are actually isometrically equivalent, meaning that there exists an isometric diffeomorphism between the two manifolds. 

Despite the popularity of these two Finsler models, this equivalence is not well established in the literature. In fact, we found only one paper referring to the isometry map from the Finslerian Poincar\'e disk onto the Funk model in the context of Zermelo's navigation problem, using polar coordinates, see Bao and Robles \cite[p. 240]{BR}. 

Therefore, the first objective of the paper is to describe in more detail the isometrical equivalence of the models (\textbf{F}) and (\textbf{P}). 
Next, we introduce a new $2$-dimensional analytic Randers model, namely
\begin{enumerate}
	\item[(\textbf{H}):]\label{H} the \textit{Finsler-Poincar\'e upper half plane}, which turns out to be precisely the Randers-type perturbation of the standard \textit{hyperbolic upper half plane}, see Loustau \cite[Section 8.2]{L20} or Stahl \cite[Chapter 4]{St93}.
\end{enumerate}	

Note that e.g. Rutz and McCarthy \cite{RM93} also considered a small perturbation of the Riemannian upper half plane, nevertheless, the metric obtained was not equivalent with the Finsler structures (\textbf{F}) and (\textbf{P}).     

In our case however, as a main result, we are able to prove that the three Finsler models (\textbf{F}), (\textbf{P}) and (\textbf{H}) are \textit{all isometrically equivalent}. This phenomena is in concordance with the behavior of the hyperbolic model spaces, as the  Riemannian  counterparts of these three models are also isometric manifolds, see e.g. Cannon, Floyd, Kenyon and Parry \cite{CFKP97}.

The isometry of the three Randers spaces reveals many interesting consequences. Most importantly, it implies that all the metric related properties which are enjoyed by one particular model can be easily proved to hold on the other two manifolds as well. In particular, based on Krist\'aly \cite{Kristaly}, we find that the first Dirichlet eigenvalue $\lambda_F$ associated to the Finsler-Laplace operator $-\Delta_F$ is zero in the case of both Finsler-Poincar\'e models (\textbf{P}) and (\textbf{H}). This provides new examples of simply connected, non-compact Finsler manifolds with constant negative flag curvature having zero first eigenvalue, which is an unexpected result considering its Riemannian counterpart proven by McKean \cite{McKean}.

The organization of the paper is the following. The next section provides a brief review of the notions of Finsler geometry used to establish our results. Section \ref{models} presents in detail the three Randers models in question. Section \ref{proof} contains the proof that the spaces (\textbf{F}), (\textbf{P}) and (\textbf{H}) are isometric. Finally, section \ref{corollaries} provides an interesting application of the results obtained.

\section{Preliminaries}   \label{preliminaries}

In this section we recall the basic notions of Finsler manifolds and Randers spaces, for further details see e.g. Bao, Chern and Shen \cite{BCS},  Ohta and Sturm \cite{OS}, and Shen \cite{Shen01}.

Let $M$ be an $n$-dimensional differentiable manifold. 
The tangent bundle of $M$ is the collection of all vectors tangent to $M$, i.e.
$$TM=\cup_{x \in M}\{(x,v): v \in T_{x} M\},$$ 
where $T_{x} M$ denotes the tangent space to $M$ at the point $x$.

The function $F: TM \to [0,\infty)$ is called a Finsler metric if it satisfies the following conditions:
\begin{enumerate}[label=(\roman*)]
	\item $F \in C^{\infty}(TM \setminus \{ 0 \})$; 
	\item $F(x,\lambda v) = \lambda F(x,v)$, for all $\lambda \geq 0$ and $(x,v) \in TM$;
	\item the Hessian matrix $\left[\left(  \frac{1}{2}F^{2}(x,v)\right) _{v^{i}v^{j}} \right]_{i,j=\overline{1,n}}$
	is positive definite for every $(x,v)\in TM \setminus \{0\}.$
\end{enumerate}
In this case we say that $(M,F)$ is a Finsler manifold.

If, in addition, $F(x,\lambda v) = |\lambda| F(x,v)$ holds for all $\lambda \in \mathbb{R}$ and $(x,v) \in TM$, then the Finsler manifold is called reversible. Otherwise, $(M,F)$ is said to be nonreversible.

The co-Finsler metric $F^*:T^*M \to [0,\infty)$ is defined as the dual metric of $F$, i.e.
\begin{equation*}  
F^*(x,\alpha) = \sup_{v \in T_xM \setminus \{0\}} ~ \frac{\alpha(v)}{F(x,v)}, \quad \forall (x,\alpha) \in T^*M,
\end{equation*}
where $T^*M = \bigcup_{x \in M}T^*_{x} M $ is the cotangent bundle of $M$ and $T^*_{x} M$ is the dual space of $T_{x} M$.

In local coordinates, the Legendre transform $J^*:T^*M \to TM$ is defined by
\begin{equation*}  
J^*(x,\alpha) = \sum_{i=1}^n \frac{\partial}{\partial \alpha_i}\left(\frac{1}{2} F^{*2}(x,\alpha)\right)\frac{\partial}{\partial x^i}.
\end{equation*}
In particular, $F(J^*(x,\alpha)) = F^*(x,\alpha)$.

If $u \in C^1(M)$, the gradient of $u$ is defined as 
\begin{equation*}  
\nabla_F u(x) = J^*(x, Du(x)), \ \forall x \in M,
\end{equation*}
where $Du(x) \in T_x^*M$ denotes the differential of $u$ at the point $x$. Note that in general, $\nabla_F$ is nonlinear. 

Given $u \in C^2(M)$, the Finsler-Laplace operator  
$\Delta_F$ is given by 
$$\Delta_F u = \mathrm{div}_F(\nabla_F u),$$
where 
\begin{equation*} 
\mathrm{div}_F(V) = 
\frac{1}{\sigma_F(x)} \sum_{i=1}^n \frac{\partial}{\partial x^i} \Big(\sigma_F(x) V^i \Big)
\end{equation*}
for some vector field $V$ on $M$, and $\sigma_F(x)$ is the density function defined by
$\sigma_F(x) = \frac{ \omega_n}{\mathrm{Vol}(B_x(1))}.$
Here $\omega_n$ and $\mathrm{Vol}(B_x(1))$ denote the Euclidean volume of the $n$-dimensional unit ball and the set 
$$B_x(1) = \Big\{ (v^i) \in \mathbb{R}^n :~ F \Big(x, \sum_{i=1}^n v^i \frac{\partial}{\partial x^i} \Big) < 1 \Big\} \subset \mathbb{R}^n,$$
respectively.
Again, the Finsler-Laplace operator $\Delta_F$ is usually nonlinear.

The Busemann-Hausdorff volume form is defined as  
\begin{equation*} \label{Hausdorff_measure}
dv_F(x) = \sigma_F(x) d x^1 \land \dots \land d x^n .
\end{equation*}

The operators $\mathrm{div}_F$ and $\Delta_F$ can be defined in a distributional sense as well, see Ohta and Sturm \cite{OS}. E.g. for every $u \in H_{\mathrm{loc}}^1(M)$, $\Delta_F u$ is defined in the weak sense as
\begin{equation*}  \label{divergence_theorem}
\int_{M} v \Delta_F u ~ dv_F(x) = -\int_{M} Dv (\nabla_F u) dv_F(x), 
\end{equation*}
for all $v \in C^{\infty}_0(M)$.

Now, if $g$ is a Riemannian metric on $M$ and the Finsler structure $F: TM \to [0,\infty)$ is given by the specific form
\begin{equation*}
F(x,v) = \sqrt{g_x(v,v)} + \beta_x(v), \quad \forall (x,v) \in TM,
\end{equation*}
where, for every $x \in M$, $\beta_x$ is a $1$-form on $M$ such that
\begin{equation}\label{def:Randers}
|\beta_x|_g = \sqrt{g^*_x(\beta_x,\beta_x)} < 1,
\end{equation}
then $F$ is called a Randers metric and $(M,F)$ is a Randers space. Here, the co-metric $g_x^*$ can be identified by the inverse of the symmetric, positive definite matrix $g_x$, induced by the Riemannian metric $g$.

Clearly, the Randers
space $(M,F)$ is reversible if and only if $\beta=0$, i.e. $(M,F)=(M,g)$ is the original Riemannian manifold. 

Finally, given two Finsler manifolds $(M_1, F_1)$ and $(M_2, F_2)$, we say that $f: M_1 \to M_2$ is an isometry if $f$ is a diffeomorphism and 
$$F_1(x,v) = F_2(f(x), Df_x(v)), \ \forall (x,v) \in TM_1,$$
where $Df_x$ denotes the differential of $f$ at the point $x$.

\section{Three models of Randers spaces} \label{models}
In this section we specify the metrics of three analytic Finslerian models of Randers type, namely the Funk model, the Finsler-Poincar\'e disk and the Finsler-Poincar\'e upper half plane. For simplicity of presentation, we consider the $2$-dimensional versions of the Randers spaces in question. 

In the sequel we use the following notations:
\begin{itemize}
	\item $D = \{(x_1, x_2) \in \mathbb{R}^2 : x_1^2+x_2^2 < 1\}$ is the $2$-dimensional Euclidean open unit disk;
	\item $H = \{(x_1, x_2) \in \mathbb{R}^2 : x_2 > 0\}$ denotes the Euclidean upper half plane;
	\item $|\cdot|$ and $\langle \cdot,\cdot\rangle$ denote the standard Euclidean norm and inner product on $\mathbb{R}^2$. 
\end{itemize}

\subsection{The Finslerian Funk model \textup{(\textbf{F})}}  
The Finslerian Funk metric $F_F :D \times \mathbb{R}^2 \to \mathbb{R}$ is given by
\begin{equation}\label{Funk}
F_F(x,v)=\frac{\sqrt{ (1-|x|^2)|v|^2 + \langle x,v\rangle^2}}{1-|x|^2}+\frac{\langle x,v\rangle}{1-|x|^2}, 
\end{equation}
for all $(x,v) \in TD$.
The pair $(D, F_F)$ is called the Funk model, which is a non-reversible Randers space having constant negative flag curvature $-\frac{1}{4}$, see Shen \cite[Example 1.3.4 \& Example 9.2.1]{Shen01}, and Cheng and Shen \cite[Example 2.1.2]{CS12}. Note that if we ommit the $1$-form $\frac{\langle x,v\rangle}{1-|x|^2}$ in \eqref{Funk}, we recover the Riemannian Klein metric, which appears in the well-known  Beltrami-Klein model having constant sectional curvature $-1$, see Loustau \cite[Section 6.2]{L20}.

\subsection{The Finsler-Poincar\'e disk \textup{(\textbf{P})}}
The Finsler-Poincar\'e metric on the open disk $D$
is defined as $F_P :D \times \mathbb{R}^2 \to \mathbb{R}$, 
\begin{equation}\label{Poincare}
F_P(x,v)=\frac{2|v|}{1-|x|^2} + \frac{4\langle x,v\rangle}{1-|x|^4}, 
\end{equation}
for every pair $(x,v) \in TD$. The Randers space $(D, F_P)$ is the famous Finsler-Poincar\'e model investigated by Bao, Chern and Shen \cite[Section 12.6]{BCS}. Again, by omitting the second term of \eqref{Poincare}, the metric reduces to the usual Riemannian Poincar\'e model, which is another well-known hyperbolic manifold of constant sectional curvature $-1$, see Loustau \cite[Section 8.1]{L20}.

\subsection{The Finsler-Poincar\'e upper half plane \textup{(\textbf{H})}}
Let us define the Finsler-Poincar\'e upper half plane model by the pair $(H, F_H)$, where $H$ is the Euclidean upper half plane and $F_H :H \times \mathbb{R}^2 \to \mathbb{R}$ is given by
\begin{equation}\label{Half}
F_H(x,v)=\frac{|v|}{x_2} + \frac{\langle w(x),v\rangle}{x_2(4+|x|^2)}, 
\end{equation}
where $w(x) \coloneqq (2x_1x_2, x_2^2-x_1^2-4)$, for all $x = (x_1, x_2) \in H$. Note that the first term in \eqref{Half} is actually the Lobachevsky metric, see Loustau \cite[Section 8.2]{L20}. Thus $F_H$ turns out to be a Randers-type perturbation of the Riemannian Poincar\'e upper half plane, another standard model of the $2$-dimensional hyperbolic space, having sectional curvature $-1$. 

\begin{proposition}
$(H, F_H)$ is a Randers space.
\begin{proof}
	It is enough to show that $|\beta_H(x)|_{g_h} <1$, where
	$$	\beta_H(x) = \frac{1}{x_2(4+|x|^2)} w(x), \quad \text{for all } x = (x_1, x_2) \in H$$ 
	and $g_h$ denotes the Riemannian metric of the Lobachevsky upper half plane, see doCarmo \cite[p. 73]{doCarmo}. 
	
	Using definition \eqref{def:Randers}, we obtain that 
	$$|\beta_H(x)|_{g_h} = \frac{|w(x)|}{4+|x|^2} < 1 , \quad \forall x \in H.$$ 
\end{proof}
\end{proposition}

\section{Main results} \label{proof}

\subsection{Equivalence of models \textup{(\textbf{P})} and \textup{(\textbf{F})} }

\begin{theorem}\label{isometry1}
Let us consider the diffeomorphism
\begin{equation*}
f: D \to D, \ f(x) = \frac{2x}{1+|x|^2},
\end{equation*}
and its inverse 
\begin{equation*}
f^{-1}: D \to D, \ f^{-1}(x) = \frac{x}{1+\sqrt{1-|x|^2}}.
\end{equation*}
Then $f$ is an isometry between the Finsler-Poincar\'e disk $(D,F_P)$ and the Funk model $(D, F_F)$. 

\begin{proof}
	It is enough to prove that
	\begin{equation}
	F_P(x,v) = F_F(f(x), Df_x(v)), \ \forall (x,v) \in TD,
	\end{equation}
	where $Df_x$ denotes the differential of $f$ at the point $x$.
	
	Given a point $x = (x_1,x_2) \in D$, the differential function $Df_x$  is determined by the Jacobian 
	\begin{equation*}
	\textbf{J}f(x) = 
	\frac{2}{(1+|x|^2)^2}
	\begin{bmatrix}
	1+|x|^2-2x_1^2 & -2x_1x_2 \\
	-2x_1x_2 & 1+|x|^2-2x_2^2
	\end{bmatrix}.
	\end{equation*}
	Then for every $v \in T_xD \cong \mathbb{R}^2$ we have 
	\begin{equation*}
	Df_x(v) = 
	\frac{2}{(1+|x|^2)^2}
	\begin{bmatrix}
	v_1(1+|x|^2)-2x_1\langle x, v\rangle \\
	v_2(1+|x|^2)-2x_2\langle x, v\rangle
	\end{bmatrix}.
	\end{equation*}
	
	Let us denote by 
	\begin{equation}\label{alfa1}
	\alpha_F(x,v) = \frac{\sqrt{ (1-|x|^2)|v|^2 + \langle x,v\rangle^2}}{1-|x|^2}
	\end{equation}
	and 
	\begin{equation}\label{beta1}
	\beta_F(x,v) = \frac{\langle x,v\rangle}{1-|x|^2}
	\end{equation}
	the norm induced by the Klein metric and the $1$-form of the Funk metric \eqref{Funk}, respectively.
	
	Expressing the terms 
	\begin{align*}
	1-|f(x)|^2 &= \frac{(1-|x|^2)^2}{(1+|x|^2)^2}, 	\\
	|Df_x(v)|^2 &= \frac{4}{(1+|x|^2)^4} \left[(1+|x|^2)^2 |v|^2-4\langle x,v\rangle^2 \right],\\
	\langle f(x),Df_x(v)\rangle &= 4\frac{1-|x|^2}{(1+|x|^2)^3} \langle x,v\rangle
	\end{align*} 
	separately, then substituting into \eqref{alfa1} and \eqref{beta1} yields
	\begin{equation*}
	\alpha_F(f(x), Df_x(v)) = \frac{2|v|}{1-|x|^2} 
	\end{equation*}
	and
	\begin{equation*}
	\beta_F(f(x), Df_x(v)) = \frac{4 \langle x,v\rangle}{1-|x|^4},
	\end{equation*}
	which concludes the proof.
\end{proof}
	
\end{theorem}

\subsection{Equivalence of models \textup{(\textbf{F})} and \textup{(\textbf{H})} }

\begin{theorem}\label{isometry2}
Let us consider the diffeomorphism
\begin{equation*}
g: D \to H, \ g(x) = \left(\frac{2x_2}{1+x_1}, \frac{2\sqrt{1-|x|^2}}{1+x_1} \right)
\end{equation*}
with its inverse function
\begin{equation*}
g^{-1}: H \to D, \ g^{-1}(x) = \left(\frac{4-|x|^2}{4+|x|^2}, \frac{4x_1}{4+|x|^2} \right).
\end{equation*}

Then $g$ is an isometry between the Funk model $(D, F_F)$ and the Finsler-Poincar\'e upper half plane $(H,F_H)$. 

\begin{proof}
	We prove that
	\begin{equation}
	F_F(x,v) = F_H(g(x), Dg_x(v)), \ \forall (x,v) \in TD.
	\end{equation}

	The Jacobian matrix of $g$ is given by 
	\begin{equation*}
	\textbf{J}g(x) = 
	-\frac{2}{(1+x_1)^2 }
	\begin{bmatrix}
	x_2 & -(1+x_1) \\
	\frac{x_1-x_2^2+1}{\sqrt{1-|x|^2}} & \frac{x_2(1+x_1)}{\sqrt{1-|x|^2}}
	\end{bmatrix}.
	\end{equation*}
	
	The Riemannian term and the $1$-form of the Finsler-Poincar\'e metric \eqref{Half} on the upper half plane $H$ is defined by
	\begin{equation*}%\label{alfa_beta2}
	\alpha_H(x,v) = \frac{|v|}{x_2} \ \text{ and } \
	\beta_H(x,v) = \frac{\langle  w(x),v\rangle}{x_2(4+|x|^2)}, 
	\end{equation*}
	where $w(x) = (2x_1x_2, x_2^2-x_1^2-4)$, for all $x = (x_1, x_2) \in H$.
	
	Expressing the term 
	\begin{equation*}
	|Dg_x(v)|^2 = 4 \frac{(1-|x|^2) |v|^2 + \langle x,v\rangle^2 }{(1+x_1)^2 (1-|x|^2)} ,	
	\end{equation*} 
	it follows that 
	\begin{align*}
	\alpha_H(g(x), Dg_x(v)) &= \frac{1+x_1}{2\sqrt{1-|x|^2}} \cdot \frac{2 \sqrt{(1-|x|^2) |v|^2 + \langle x,v\rangle^2} }{(1+x_1)\sqrt{1-|x|^2}} \\
	& = \frac{\sqrt{ (1-|x|^2)|v|^2 + \langle x,v\rangle^2}}{1-|x|^2} \\
	& = \alpha_F(x,v),
	\end{align*}
	while for the $1$-form $\beta_H$ we use the following calculations:
	\begin{align*}
	w(g(x)) &= \left( 8\frac{x_2 \sqrt{1-|x|^2}}{(1+x_1)^2}, -8 \frac{|x|^2 + x_1}{(1+x_1)^2} \right), \\
	4+|g(x)|^2 &= \frac{8}{1+x_1}.
	\end{align*}		 
	After a direct computation we obtain that
	\begin{equation*}
	\langle w(g(x)), Dg_x(v) \rangle = \frac{16 \langle x, v \rangle }{(1+x_1)^2 \sqrt{1-|x|^2}}, 
	\end{equation*}	 
	thus
	\begin{align*}
	\beta_H(g(x), Dg_x(v)) &= \frac{(1+x_1)^2}{16\sqrt{1-|x|^2}} \cdot \frac{16 \langle x, v \rangle }{(1+x_1)^2 \sqrt{1-|x|^2}} \\
	& = \frac{\langle x,v\rangle}{1-|x|^2} \ 
	= \ \beta_F(x, v) .
	\end{align*}		
\end{proof}
	
\end{theorem}

\subsection{Equivalence of models \textup{(\textbf{H})} and \textup{(\textbf{P})} }

\begin{theorem}\label{isometry3}
Let us consider the diffeomorphism
\begin{equation*}
h: H \to D, \ h(x) = \left(\frac{4-|x|^2}{|x|^2+4 x_2+4}, \frac{4x_1}{|x|^2+4x_2+4} \right),
\end{equation*}
and its inverse 
\begin{equation*}
h^{-1}: D \to H, h^{-1}(x) = \left(\frac{4x_2}{|x|^2+2x_1 + 1}, \frac{2-2|x|^ 2}{|x|^2+2 x_1 + 1} \right).
\end{equation*}
Then $h$ is an isometry between the Finslerian upper half plane $(H, F_H)$ and the Finsler-Poincar\'e disk $(D, F_P)$. 
	
	\begin{proof}
	It is enough to show that
	\begin{equation}
	F_H(x,v) = F_P(h(x), Dh_x(v)), \ \forall (x,v) \in TH.
	\end{equation}
		
	The Jacobian of $h$ can be written as 
	\begin{equation*}
	\textbf{J}h(x) = 
	\frac{-4}{(|x|^2+4x_2+4)^2 }
	\begin{bmatrix}
	2x_1(x_2+2) & (x_2+2)^2-x_1^2 \\
 	x_1^2 - (x_2+2)^2 & 2x_1(x_2+2)
	\end{bmatrix}_.
	\end{equation*}
	
	Let us denote by
	\begin{equation*}%\label{alfa_beta3}
	\alpha_P(x,v) = \frac{2|v|}{1-|x|^2} \ \text{ and } \
	\beta_P(x,v) = \frac{4\langle x,v\rangle}{1-|x|^4} 
	\end{equation*}
	the terms determined by the inner product and the $1$-form of the Finsler-Poincar\'e metric \eqref{Poincare} on the disk $D$.
	
	First we compute the following terms:
	\begin{align*}
	1-|h(x)|^2 &= \frac{8 x_2}{|x|^2+4x_2+4}, 	\\
	1+|h(x)|^2 &= 2\frac{|x|^2+4}{|x|^2+4x_2+4}, 	\\
	1-|h(x)|^4 &= \frac{16 x_2(|x|^2+4)}{(|x|^2+4x_2+4)^2}, 	\\
	|Dh_x(v)| &= \frac{4|v|}{|x|^2+4x_2+4}, \\ 
	\langle h(x), Dh_x(v) \rangle &= \frac{-4}{(|x|^2+4x_2+4)^3} \cdot \\
	&\Big\{ (x_1^2-(x_2+2)^2)\big(4x_1v_1 - (4-|x|^2)v_2\big) \\ 
	&+ 2x_1(x_2+2)\big((4-|x|^2)v_1+4x_1v_2\big)\Big\} \\
	& = 4 \cdot\frac{2x_1x_2v_1 + (x_2^2-x_1^2-4)v_2}{(|x|^2+4x_2+4)^2} ,
	\end{align*}  
	for every $v=(v_1, v_2) \in T_xH \cong \mathbb{R}^2$.	
	It follows that 
	\begin{equation*}
	\alpha_P(h(x), Dh_x(v)) = 
	\frac{2|Dh_x(v)|}{1-|h(x)|^2} = \frac{|v|}{x_2}
	\end{equation*}
	and 
	\begin{align*}
	\beta_P(h(x), Dh_x(v)) &= \frac{2x_1x_2v_1 + (x_2^2-x_1^2-4)v_2}{x_2(|x|^2+4)} \\
	&= \beta_H(x, v),
	\end{align*}		
	which concludes the proof.
\end{proof}
\end{theorem}

\begin{remark}
	Note that for the previous isometries we have 
	$$h^{-1}=g \circ f,$$ 
	i.e. the following diagram is commutative:
	\begin{center}
	\includegraphics[width=0.22\textwidth]{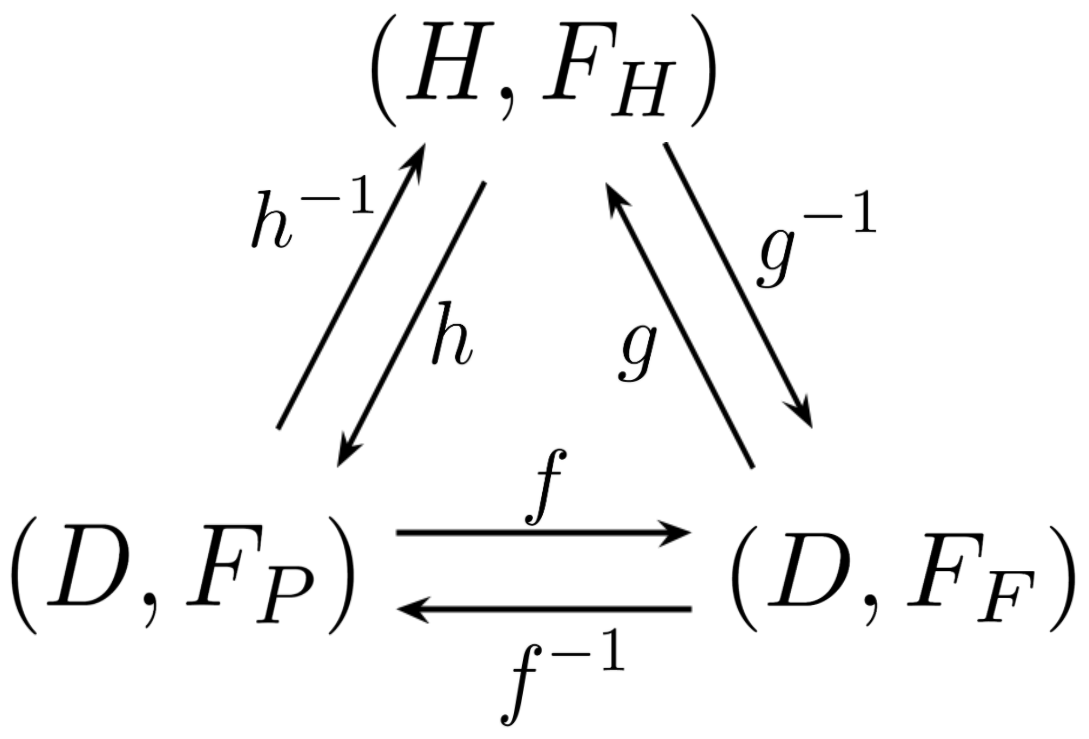}
	\end{center}
	
	Moreover, these diffeomorphisms actually coincide with the appropriate isometries available between the Riemannian counterpart of the models, i.e. the Beltrami-Klein disk, the Riemannian Poincar\'e disk, and the hyperbolic upper half plane, see e.g. Cannon, Floyd, Kenyon and Parry \cite[p. 69]{CFKP97}.  
	This is illustrated by the proofs of Theorems \ref{isometry1}--\ref{isometry3}, as well, where the norms $\alpha_P, \alpha_F, \alpha_H$ and the $1$-forms $\beta_P, \beta_F, \beta_H$ turn out to be the pullbacks of one another by the corresponding isometries $f,g$ and $h$. 
\end{remark}

\section{Consequences} \label{corollaries}

An important byproduct of the isometries given in Theorems \ref{isometry1}, \ref{isometry2} and \ref{isometry3} is the fact that all the metric related properties of one of the models can be easily transferred to the other two manifolds by the appropriate isometry function. 
 
To give an interesting example, let us consider the first eigenvalue associated to the Finsler-Laplace operator $-\Delta_F$ on the spaces $\textup{(\textbf{F})}$, $\textup{(\textbf{P})}$ and $\textup{(\textbf{H})}$, respectively. 

Given a Finsler manifold $(M,F)$, the first eigenvalue associated to $-\Delta_F$ (also called the fundamental frequency) is defined as
$$\lambda_{1,F}(M) = \inf_{u \in H_{0,F}^1(M) \setminus \{0\}} \frac{\int_M {F^*}^2(x, Du(x)) \text{d}v_F(x)}{\int_M u^2(x) \text{d}v_F(x)} ,$$
where $H_{0,F}^1(M)$ is the closure of $C_0^\infty(M)$ with respect to the norm 
$$\|u\|_{ \scriptscriptstyle{H_{0,F}^1}} = \left(\int_M {F^*}^2(x, Du(x)) \text{d}v_F(x) + \int_M u^2(x) \text{d}v_F(x) \right)_,^{\frac{1}{2}} $$
see Ge and Shen \cite{GeShen01}, Ohta and Sturm \cite{OS}.

According to Krist\'aly \cite[Theorem 1.3]{Kristaly}, in case of the Finslerian Funk model $(D,F_F)$, we have that
$$\lambda_{1,F_F}(D) = 0.$$

Combining this with the isometries proven in Theorems \ref{isometry1} and \ref{isometry2}, we obtain the following result:

\begin{corollary}\label{Corollary}
In case of the Finsler-Poincar\'e disk $(D,F_P)$ and the Finslerian upper half plane $(H, F_H)$, we have 
$$\lambda_{1,F_P}(D) = \lambda_{1,F_H}(H) = 0.$$
\end{corollary}  

These assertions are in sharp contrast with the result of McKean \cite{McKean}, which states that for every complete, $n$-dimensional, simply connected Riemannian manifold $(M,g)$ having sectional curvature bounded above by $-\kappa^2 (\kappa>0)$, one has the following spectral gap:
$$\lambda_{1,g}(M) \geq \frac{(n-1)^2}{4}\kappa^2.$$
In fact, on the Beltrami-Klein disk and the Riemannian upper half plane the first eigenvalue is precisely $\frac{1}{4}$, since in the case of the $n$-dimensional hyperbolic space $(\mathbb{H}^n, g_h)$ of constant curvature $-\kappa^2 (\kappa>0)$, we have 
$$\lambda_{1,g_h}(\mathbb{H}^n) = \frac{(n-1)^2}{4}\kappa^2,$$
see Chavel \cite[p. 46]{Chavel}.

Therefore,  Corollary \ref{Corollary} provides another example highlighting the anisotropic nature of the Finsler metrics $F_P$ and $F_H$, despite the simplicity of these models.

\section*{Acknowledgment} 
The authors are supported by the National Research, Deve\-lopment and Innovation Fund of Hungary, financed under the K$\_$18 funding scheme, Project No.  127926.

\end{document}